\documentclass[11pt]{article}

\usepackage{amsfonts,amssymb}
\usepackage[all, knot]{xy}

\newcommand{\Z}{\mathbb{Z}}

\newtheorem{theorem}{Theorem}[section]

\newtheorem{lemma}[theorem]{Lemma}
\newtheorem{proposition}[theorem]{Proposition}

\newtheorem{example}[theorem]{Example}

\newtheorem{remark}[theorem]{Remark}

\input{epsf.sty}

\begin{document}

\title{Polynomial Cocycles of Alexander Quandles and Applications}

\author{Kheira Ameur \quad and \quad
Masahico Saito\footnote{Supported in part by NSF Grant DMS \#0301089,
 \#0603876.}
\\ University of South Florida
}

\maketitle

\begin{abstract}
Cocycles are constructed  by polynomial expressions
 for Alexander quandles. As applications, non-triviality of
some quandle homology groups are proved, and quandle cocycle invariants
of knots are studied.
In particular, for 
an infinite family of quandles,
the non-triviality of quandle homology groups is proved for all odd dimensions.
\end{abstract}

\section{Introduction}

Quandles are self-distributive sets with additional
properties (see below for details). They  have been used in
the study of knots since 1980s. Cohomology theories
 of quandles  were developed as  modifications
of rack cohomology theory~\cite{FRS}, and their cocycles have
been used to construct
invariants of knots and knotted surfaces~\cite{CJKLS}.
Quandle cohomology theory was further generalized
\cite{CEGS,CES,FRS}.

To compute the quandle cocycle invariant, explicit  cocycles are
required. At first the computations relied on cocycles found by
computer calculations. A significant progress was made in
computations of the invariants after Mochizuki discovered a family
of 2- and 3-cocycles for dihedral and other linear Alexander
quandles, for which cocycles were written by polynomial
expressions. Formulas for important families of knots and knotted
surfaces and their applications followed \cite{AS,Hatake,Iwakiri}.
Homology groups for higher dimensions are studied in \cite{NP} for
dihedral quandles.

In this paper, following the method of the construction by
Mochizuki, a variety of $n$-cocycles for $n\geq 2$ are constructed
for some Alexander quandles  by polynomial expressions.
As  applications,
we prove the
non-triviality of
 quandle homology groups
 for an infinite family of Alexander quandles for
all odd dimensions.
These cocycles are also used  to compute the invariants for
some families of knots.
Much of the results are based on the Ph.D. dissertation by K.A. \cite{KAdissertation}.

The paper is organized as follows.
In Section~\ref{prelimsec}, a brief review of quandle homology groups
are given and a key lemma is proved on polynomial cocycles.
Quandle cocycle invariants of knots are studied in Section~\ref{knotsec}, and
non-triviality of quandle homology groups is proved in Section~\ref{homgpsec}.

\section{Preliminaries}\label{prelimsec}

A {\it quandle}, $X$, is a set with a binary operation
$(a, b) \mapsto a * b$
satisfying the three conditions:
(I) For any $a \in X$,
$a* a =a$,
(II) for any $a,b \in X$, there is a unique $c \in X$ such that
$a= c*b$, and
(III)
for any $a,b,c \in X$, we have
$ (a*b)*c=(a*c)*(b*c). $
A {\it rack} is a set with a binary operation that satisfies
(II) and (III).
Racks and quandles have been studied in, for example,
\cite{Br88,FR,Joyce,Matveev}.
For  $\Lambda ={\Z }[t, t^{-1}]$,
any $\Lambda$-module $M$
is a quandle with
$a*b=ta+(1-t)b$, $a,b \in M$,
that is
called an {\it  Alexander  quandle}.

\begin{sloppypar} 
 A cohomology theory
 of quandles  was defined~\cite{CJKLS} as  a modification
of rack cohomology theory~\cite{FRS} as follows.
Let $C_{n}^R(X)$ be the free abelian group generated by $n$-tuples
$(x_1,\ldots,x_n)$ of elements of a quandle $X$. Define a
homomorphism $\partial_n: C_{n}^R(X) \longrightarrow C_{n-1}^R(X)$
by
\begin{eqnarray*}
\lefteqn{\partial_{n}(x_1,\ldots,x_n)=\sum_{i=2}^{n}(-1)^i[(x_1,x_2,\dots,x_{i-1},x_{i+1},\ldots,x_{n})}\\
&-&(x_1*x_i,x_2*x_i,\ldots,x_{i-1}*x_i,x_{i+1},\ldots,x_{n})]
\end{eqnarray*}
for $n\geq 2$ and $\partial_{n}=0$ for $n\leq 1$. Then
$C_{*}^R(X)=\{C_{n}^R(X),\partial_n\}$ is a chain complex. Let
$C_{n}^D(X)$ be the subset of $C_{n}^R(X)$ generated by n-tuples
$(x_1,\ldots, x_n )$ with $x_i=x_{i+1}$ for some $i \in
\{1,\ldots,n-1\}$ if $n\geq 2$; otherwise let $C_{n}^D=0.$ If $X$ is
a quandle, then $\partial_n(C_{n}^D)\subset\partial_{n-1}(C_{n}^D)$
and $C_{*}^D(X)=\{C_{n}^D(X),\partial_n\}$ is a subcomplex of $
C_{*}^R(X)$. Put $ C_{n}^{Q}=C_{n}^R(X)/ C_{n}^D(X)$ and
$C_{*}^Q(X)=\{C_{n}^Q(X),\partial'_n\}$ where, $\partial'_n$ is the
induced homomorphism. Henceforth, all boundary maps may be denoted
by $\partial_{n}$. The superscripts $R,Q$ and $D$, respectively,
represent rack, quandle, and degenerate chain complexes. For an
Abelian group $G$, define the chain and the cochain
complexes by 
 \begin{eqnarray*}
C_{*}^W(X;G)=C_{*}^W(X)\otimes G  &  & \partial=\partial
\otimes {\rm id}, \\
C_{W}^*(X;G)={\rm Hom}(C_{*}^W(X),G) & & \delta={\rm Hom}(\partial,id)
\end{eqnarray*}
in the usual way, where $W=D,R,Q.$ 
The groups of cycles and
boundaries are denoted respectively by
$ker(\partial)=Z_{n}^W(X;G)\subset C_n^W(X;G)$ and ${\rm
Im}(\partial)=B_W^n(X;G)\subset C_W^n(X;G)$ while the cocycles and
coboundaries are denoted respectively by  $Z_W^n(X;G)$ and $B_W^n(X;G)$, respectively.
The $n\/$th quandle homology
group with coefficient group $G$ is defined by
$$ H_{n}^{Q}(X;G)=H_{n}(C_{*}^w(X;G))=Z_{n}^Q(X;G)/B_n^Q(X;G),$$
and the 
quandle cohomology group with coefficient group $G$,
$H_{Q}^{n}(X;G)$,  is defined similarly.
The following is the Key Lemma of the paper.
\end{sloppypar} 

\begin{lemma} \label{polycocylemma}
\begin{sloppypar}
Let $X=\Z_p[t, t^{-1} ]/ (g(t))$ for some prime $p$.
 Let
$a_i=p^{m_i}$, for $i=1,\ldots,n-1$, where $m_i$
are non-negative integers. For a positive integer $n$,
let  $f:X^n \rightarrow X$ be defined by
$$f(x_{1},x_{2},\ldots,x_{n})
=(x_{1}-x_{2})^{a_1}(x_{2}-x_{3})^{a_2}\cdots(x_{n-1}-x_{n})^{a_{n-1}}x_n^{a_n}.$$
\begin{enumerate}
\setlength{\itemsep}{-3pt}
\item If $a_n=0$, then $f$ is an $n$-cocycle $(\in Z^n_{\rm Q}(X;X))$.
\item  If $a_n=p^{m_n}$ for a positive integer ${m_n}$, then $f$ is an
$n$-cocycle if  $g(t)$ divides $1-t^a$, where
$a=a_1+a_2+\cdots+a_{n-1}+a_n$.
\end{enumerate}
\end{sloppypar}
\end{lemma}
{\it Proof.\/} From the definition of $\delta$,  for $i=1, \ldots ,
n$,  
 we compute using  the notation $y_i=x_i - x_{i+1}$ 
\begin{eqnarray*}
 \lefteqn{\delta  f(x_1,\ldots ,x_{n+1})}\\
 &=& \sum^{n+1}_{i=2}(-1)^{i}[f(x_1,x_2,\ldots ,x_{i-1},x_{i+1}, \ldots, x_{n+1})
 -  f(x_1*x_i,\ldots,x_{i-1}*x_{i},x_{i+1},\ldots , x_{n+1}) ]  \\
 &=& \sum^{n}_{i=2}(-1)^{i}y^{a_1}_1y^{a_2}_2\cdots
 (y_{i-1}+y_i)^{a_{i-1}}y^{a_i}_{i+1}\cdots y^{a_{n-1}}_nx^{a_n}_{n+1}
 + (-1)^{n+1}y^{a_1}_1y^{a_2}_2\cdots y^{a_{n-1}}_{n-1}x^{a_n}_n\\
 & -&\sum^{n}_{i=2}(-1)^i(ty_1)^{a_1}(ty_2)^{a_2}\cdots
 (ty_{i-1}+y_i)^{a_{i-1}}\cdots y^{a_{n-1}}_nx^{a}_{n+1} \\
 & -& (-1)^{n+1}t^{a}y^{a_1}_1y^{a_2}_2\cdots y^{a_{n-1}}_{n-1}
 (t y_n + x_{n+1})^{a_n}
 \end{eqnarray*}
which simplifies to
\begin{eqnarray*}
& & (-1)^{n}y^{a_1}_1\cdots y^{a_{n-1}}_{n-1}x^{a_n}_{n+1}
- (-1)^n t^{a_1+\cdots+a_{n-1}}y^{a_1}_1\cdots
y^{a_{n-1}}_{n-1}x^{a_n}_{n+1}\\
& & + (-1)^{n+1}y^{a_1}_1 y^{a_2}_2\cdots
y^{a_{n-1}}_{n-1}x^{a_n}_n
- (-1)^{n+1}t^{a_1+\cdots a_{n-1}}y^{a_1}_1\cdots
y^{a_{n-1}}_{n-1}(ty_n+x_{n+1})^{a_n}.
\end{eqnarray*}
If  $a_n=0$ then we see $\delta  f(x_1,..,x_{n+1})=0$.
 If $a_n=p^{m_n}$ then we have
  \begin{eqnarray*}
 \lefteqn{\delta  f(x_1,\ldots ,x_{n+1})} \\
  &=&(-1)^ny^{a_1}_1\cdots
 y^{a_{n-1}}_{n-1}x^a_{n+1}+(-1)^{n+1}y^{a_1}_1y^{a_2}_2\cdots
 y^{a_{n-1}}_{n-1}x^{a_n}_n
 -(-1)^{n+1}t^{a_{1}+\cdots+a_{n-1}+a_n}y^{a_1}_1\cdots
 y^{a_n}_{n}\\
  &=&(-1)^{n+1}(1-t^{a})y^{a_1}_1\cdots
 y^{a_n}_n \quad = \quad 0 \in X
\end{eqnarray*}
by assumption. Hence $f$ is an n-cocycle. $\hfill$ $\Box$

\bigskip

For the purpose of constructing knot invariants for applications
 in later sections, we mainly use  quandle $2$- and
$3$-cocycles. The properties of these cocycles that we need for the
knot  invariants are specifically formulated as follows. A
quandle 
$2$-cocycle $\phi$ is regarded as a function $\phi: X \times X
\rightarrow A$ with the $2$-cocycle condition
$$ \phi(x,y)+\phi(x*y, z)=\phi(x,z)+\phi(x*z, y*z)$$
for all  $x,y,z \in X$, and $\phi(x,x)=0$ for all $x \in X$. A
quandle 
$3$-cocycle $\theta$ is regarded as a function $\theta: X \times X
\times X \rightarrow A$ with the $3$-cocycle condition
$$\theta(x,z,w)+\theta(x,y,z)+\theta(x*z,y*z,w)=\theta(x*z,y,w)+\theta(x,y,w)+\theta(x*w,y*w,z*w)$$
for any $x,y,z,w \in X$, and $\theta(x,x,y)=0$, $ \theta(x,y,y)=0$
for all $x,y \in X.$

\begin{figure}
\begin{center}
\mbox{
\epsfxsize=2.5in
\epsfbox{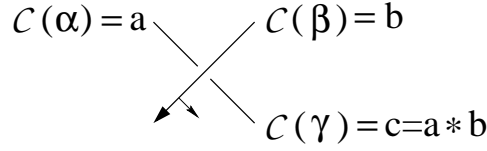}
}
\end{center}
\caption{ Quandle relation at a crossing  }
\label{qcolor}
\end{figure}

Let $X$ be a fixed quandle.
Let $K$ be a given oriented classical knot or link diagram,
and let ${\cal R}$ be the set of (over-)arcs.
The normals (normal vectors) are given in such a way that
the ordered pair
(tangent, normal) agrees with
the orientation of the plane, see Fig.~\ref{qcolor}.
A (quandle) {\it coloring} ${\cal C}$ is a map
${\cal C} : {\cal R} \rightarrow X$ such that at every crossing,
the relation depicted in Fig.~\ref{qcolor} holds.
The  (ordered) colors ${\cal C}(\alpha)$, ${\cal C}(\beta)$
are called {\it source} colors.
Let ${\rm Col}_X(K)$ denote the set of colorings of a knot diagram $K$
by a quandle $X$.

\begin{figure}[htb]
\begin{center}
\mbox{
\epsfxsize=3.2in
\epsfbox{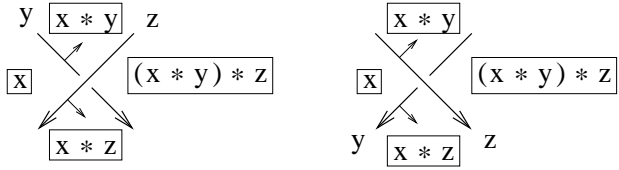}
}
\end{center}
\caption{ Quandle colorings of regions }
\label{Foxcolors}
\end{figure}

For a coloring $ {\mathcal C}$, there is a coloring of regions that
extend $ {\mathcal C}$ as depicted in Fig.~\ref{Foxcolors}. Let
$(x,y,z)=(x_{\tau}, y_{\tau}, z_{\tau})$ be the colors near a
crossing $\tau$ such that $x$ is the color of the region (called the
source region) from which both orientation normals of over- and
under-arcs point, $y$ is the color of the under-arc (called the
source under-arc) from which the normal of the over-arc points, and
$z$ is the color of the over-arc. See Fig.~\ref{Foxcolors}.

\begin{figure}[htb]
\begin{center}
\mbox{
\epsfxsize=1.5in
\epsfbox{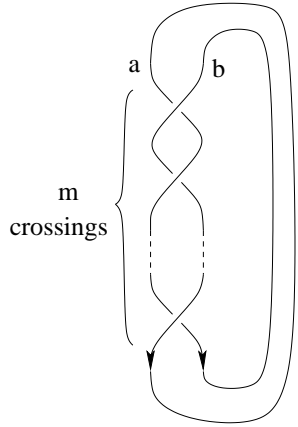}
}
\end{center}
\caption{Colorings of torus knots }
\label{torusknots}
\end{figure}

\begin{example}\label{colorex}
{\rm
Let $X$ be an Alexander quandle.
Let $T(2, m)$ be the $(2, m)$-torus knot
that is the closure of the closed braid of $2$-string braid with
$m$ positive crossings for a positive integer $m$.
For a negative integer $m$, $T(2,m)$ consists of $|m|$ negative crossings.
Denote by $\xi_k$ the polynomial
\begin{displaymath} \xi_k =\xi_k (t)=\sum^{k-1}_{i=0}(- t)^i,\end{displaymath}
and define $\xi_0=0$ as convention. Note that for $m>1$, 
the polynomial $\xi_m (t)$ 
is the Alexander polynomial $\Delta_{T(2,m)} (t)$
of the knot $T(2,m)$  (see, for example, \cite{Mura}). 

Then it is seen by induction that
 if $(a,b)\in X \times X$ is the top color vector (the elements $a$ and $b$ are
 assigned to the top left and right arcs of a $2$-string braid, respectively),
 of a coloring of
 $T(2,m)$ by $X$, then
 the
$k$\/th color vector (the pair of colors after $k$\/th crossing)
is 
$(a + \xi_k (b-a), b+ \xi_k (a-b) )$ 
where $1 \leq k
\leq m$.
In particular, any top color vector extends to a coloring
of $T(2,m)$ for the quandle
$X=\Z_p[t, t^{-1}]/ (\xi_m (t))$.

\begin{figure}[htb]
\begin{center}
\mbox{
\epsfxsize=4.5in
\epsfbox{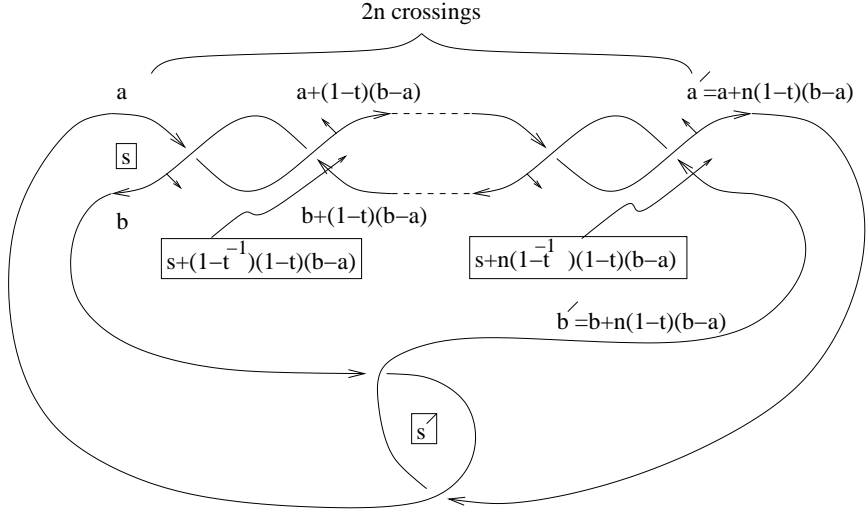}
}
\end{center}
\caption{Colorings of twist knots }
\label{twistknots}
\end{figure}

Another example, a coloring of twist knots (see for example~\cite{Rolf})
is depicted in Fig.~\ref{twistknots}.
We denote the twist knot with
 with $2n+2$ crossings
as depicted in the figure by $k(2n)$.
} \end{example}

The cocycle invariant for classical knots \cite{CJKLS}
was defined as follows. Let $\phi \in Z^2_{\rm Q}(X;A)$
be a quandle $2$-cocycle of a finite quandle $X$ with the
coefficient group $A$.
Let ${\cal C}$ be a coloring of a given knot diagram $K$ by $X$.
The {\em Boltzmann weight} $B( {\cal C}, \tau)$ at a crossing $\tau$ of $K$
is then defined by  $B( {\cal C}, \tau)=\phi(x_{\tau}, y_{\tau})^{\epsilon(\tau)}$,
where $x_{\tau}$, $y_{\tau}$ are source colors at $\tau$ and
$\epsilon(\tau)$ is the sign ($\pm 1$) of $\tau$. In Fig.~\ref{qcolor},
it is a positive crossing if the under-arc is oriented downward.
Here $B( {\cal C}, \tau)$ is an element of $A$ written multiplicatively.
The formal sum (called a state-sum) in the group ring $\Z[A]$
 defined by 
$\Phi_{\phi} (K) = \sum_{{\cal C}\in {\rm Col_X(K)}} \prod_{\tau}
B( {\cal C}, \tau ) $ is called the {\it quandle cocycle invariant}. 
The invariant is also defined 
by
$ \{ \sum_{\tau} {\epsilon(\tau)} \phi(x_{\tau}, y_{\tau}) \ | \ {\cal C}\in {\rm
Col_X(K)} \} $ as a multiset, in which case the values of the cocycles are
written additively.
If the quandle $X$ is finite, the invariant as a multiset can be
written by an expression similar to those for the state-sums
as follows. 
 Suppose 
a multiset of group elements is
given by 
$\{ \sqcup_{m_1} g_1, \ldots,
\sqcup_{m_\ell} g_\ell \} $,
where $\sqcup_{m_i} g_i$ denotes $m_i$ copies of $g_i$ 
(the positive  integer $m_i$ is called the multiplicity of $g_i$),  
then we use the polynomial notation
$m_1 U^{g_1} + \cdots + {m_\ell} U^{g_\ell }$ where
 $U$ is
 a formal symbol. For example, the multiset value of
the invariant for a trefoil with the Alexander quandle
$X=\Z_2[t, t^{-1}]/ (t^2+t+1)$ with the same coefficient $A=X$ and  a certain
$2$-cocycle is $\{ \sqcup_4 (1) ,  \sqcup_{12} (t+1) \}$, and
 is
denoted by $4 + 12 U^{(t+1)}$, where we use the convention $U^0=1$ and
exponential rules apply.

It is seen that $c= \sum_{\tau} {\epsilon(\tau)} (x_{\tau}, y_{\tau})$ is a
$2$-cycle~\cite{CKS}, and the contribution $\sum_{\tau}
{\epsilon(\tau)} \phi(x_{\tau}, y_{\tau})$ is regarded as the
evaluation $\phi(c)$ of the cycle by the cocycle $\phi$. Hence a
colored knot diagram represents a $2$-cycle, and the multiset of
evaluations over all colorings is the cocycle invariant.  The quandle cocycle
invariants have also been defined for knotted surfaces in
$4$-space, in a similar manner, using quandle $3$-cocycles
and triple points on projections.

\begin{sloppypar}
The quandle cocycle invariant using region  colorings (sometimes
called ``{\it{shadow}}'' colorings) and $3$-cocycles were considered
in \cite{FRS}.
Let $\phi \in Z^3_{\rm Q}(X;A)$ be a $3$-cocycle.
Then the weight in this case is defined by
 $B( {\cal C}, \tau)=\phi(x_{\tau}, y_{\tau}, z_{\tau})^{\epsilon(\tau)}$
 where $\epsilon(\tau)$ is $+1$ or $-1$, for a positive or a negative
crossing respectively. Then the $3$-cocycle invariant is defined
by
$\Phi_{\phi}(K) = \sum_{{\cal C}\in {\rm
Col_X(K)}} \prod_{\tau} B( {\cal C}, \tau ) $ as a state-sum, and by
$ \{ \sum_{\tau} {\epsilon(\tau)} \phi(x_{\tau}, y_{\tau}, z_{\tau})  \ | \ {\cal C}\in {\rm
Col_X(K)} \} $ as a multiset.
\end{sloppypar}

\section{Quandke Cocycle Invariants of Knots  with Polynomial Cocycles}
\label{knotsec}

In this section we exhibit examples of quandle cocycle invariants
that are obtained from polynomial cocycles.
The  following particular case  illustrates our calculations of the invariant using
polynomial cocycles, and will be used in the next section.

\begin{proposition} \label{2cocyinvprop}
Let  $X=\mathbb{Z}_p[t,t^{-1}]/ (\xi_m (t)) $ for a prime $p$ and a
positive integer $m$, and let $f : X \times X \rightarrow X$ be
defined by $f(x_1, x_2)=(x_1-x_2)^{a_1} x_2^{a_2}$ for
$a_i=p^{m_i}$, $i=1,2$, where $m_i$ are non-negative integers.
Suppose that $\xi_m (t)$ divides $1-t^{(a_1+a_2)}$ in
$\Z_p[t, t^{-1}]$. Then the cocycle invariant $ \Phi_f(T(2,m))$
of the torus knot $T(2,m)$
 is given by the multiset
$$ \Phi_f(T(2,m))=
\{ \sqcup_{|X|} 
 (t^2 \xi_{m}' (t) )^{a_2} 
s^{(a_1+a_2)} \ | \ s \in X \} ,$$
where $|X|$ denotes the number of elements of $X$,
and $\xi_m'$ is the derivative of $\xi_m$. 
If $m$ is a negative integer, then the invariant consists of the
negatives of the invariant for $T(2, |m|)$.
\end{proposition}
{\it Proof.\/} By Lemma~\ref{polycocylemma}, indeed $f \in Z^2_{\rm
Q}(X;A)$. By assumption any top color vector $(a,b)$ extends to a
coloring from Example~\ref{colorex}.
Recall also from Example~\ref{colorex} that
just below  the
$k$\/th crossing  we have the $k$\/th  color vector
$(t\xi_{k-1}+\xi_kb,t\xi_ka+\xi_{k+1}b)$, so the contribution to the
invariant is computed as
\begin{eqnarray*}
\lefteqn{\sum_{k=1}^m f(t \xi_{k-1}a +\xi_k b, t \xi_k a + \xi_{k+1} b)}\\
&=& \sum_{k=1}^m( t \xi_{k-1} a +\xi_{k}b -t\xi_{k}a -
\xi_{k+1}b)^{a_1}
(t\xi_{k}a+\xi_{k+1}b)^{a_2}\\
&=&
\sum_{k=1}^m [ (a-b)(\xi_{k+1}-\xi_{k}) ]^{a_1}[b+t\xi_{k}(a-b)]^{a_2}\\
&=& (a-b)^{a_1}b^{a_2} \sum_{k=1}^m(- t)^{k a_1} +(a-b)^{a_1+a_2}
t^{a_2} \sum_{k=1}^m(-t)^{ka_1}\xi_k^{a_2}.
\end{eqnarray*}
The first term is written as
$$ (a-b)^{a_1}b^{a_2} \left( \sum_{k=1}^m (-t)^k \right)^{a_1}
= (a-b)^{a_1}b^{a_2} [ -t \xi_m ] ^{a_1}
$$
which vanishes in $X$. Note that by
assumption $t^{a_1+a_2}=1$ in $X$, so that $t^{a_1}=t^{-a_2}$. Hence
the second term is written as
$$(a-b)^{a_1+a_2} t^{a_2} \left( \sum_{k=1}^m (-1)^k t^{-k}  \xi_k \right)^{a_2} .$$
Thus we compute $S_m=\sum_{k=1}^m (-t)^{-k} \xi_k$. By induction we obtain
$S_m=\sum_{k=0}^{m-1} (k+1)(-t)^{k-m}.$
On the other hand, we compute
\begin{eqnarray*}
S_{m}
     &=& (-t)^{-m} \left( \sum_{k=0}^{m-1}(-t)^{k}+\sum_{k=0}^{m-1}k(-t)^k \right)\\
     &=&(-t)^m \left(\xi_{m}+t\sum_{k=0}^{m-1}(-1)^kkt^{k-1} \right)\quad
     = \quad (-t)^m (\xi_{m}+t \xi_{m}').
\end{eqnarray*}
Hence the contribution of
the coloring induced  by a top color vector
$(a,b)$  is
\begin{displaymath}(a-b)^{a_1+a_2}(-t)^{a_2(2-m)}(\xi_{m}')^{a_2}=
(a-b)^{a_1+a_2}( t^2 \xi_{m}')^{a_2}, 
\end{displaymath}
since $(-t)^m=1$ in $X$.
Let $s=a-b$.
Then we have
$$\{ (a-b)^{(a_1+a_2)} \ | \ (a,b) \in X \times X \}
= \{ \sqcup_{|X|} s^{(a_1+a_2)} \ | \ s \in X  \} $$ and the result
follows.

If $m$ is negative, then all crossings are negative. Then consider
the diagram of $T(2, m)$ that is the mirror of the diagram used
above of $T(2, |m|)$, with opposite orientation. Then also consider
the colors of $T(2, m)$ at the bottom arcs $(a,b)$. Then the
contribution from the coloring induced by this bottom color vector
coincides with the negative of the original. Hence the invariant
$\Phi_f(T(2,m))$ is the multiset that  consists of the negative of
$\Phi_f(T(2,|m|))$. $\hfill$ $\Box$

\begin{example} {\rm
For
$f(x_1,x_2)=(x_1-x_2)^4 x_2$, the above formula gives
the $2$-cocycle invariant
$$\Phi_f (T(2,5))
=\{ \sqcup_{16} 0, \sqcup_{80} (t+1), \sqcup_{80} (t^3), \sqcup_{80}
(t^3+t+1) \} $$ where $m=5=2^2+1$ and
$X=\mathbb{Z}_2[t,t^{-1}]/ (\xi_5 (t) )$. } \end{example}

For $3$-cocycle invariants, similar calculations can be carried out
with the $3$-cocycle
$f(x_1,x_2,x_3)=(x_1-x_2)^{a_1}(x_2-x_3)^{a_2}$ where $a_i=
p^{m_i}$ for $i=1,2$, where
  $m_i$'s are non-negative integers.
In particular, let
$X=\Z_p[t, t^{-1}]/( \xi_m(t) )$ and  $(a,b) \in X \times X$ be a top color vector
 for  a coloring of $T(2,m)$ by $X$
 as in Proposition~\ref{2cocyinvprop},
 with the source region color $c \in X$,
 where $m$ is a positive integer.
 Then the contribution to the cocycle invariant of this coloring is
 given by
 $$  (a-b)^{a_1+a_2}\sum_{k=1}^{m}( \xi_k (t) )^{a_1} (-t)^{ka_2} . $$
 In particular, if $t^{a_1+a_2}=1$ in $X$, then the contribution is
 given by $$  (a-b)^{a_1+a_2} (-t)^{a_1(1-m)}( \xi'_m (t) )^{a_1}
 =  (a-b)^{a_1+a_2} ( - t \xi'_m (t) )^{a_1}. $$
 Thus we obtain

  \begin{proposition}\label{3cocyinvprop}
\begin{sloppypar}
 Let $p$ be a prime and $m$ be a positive integer.
 Let $X=\mathbb{Z}_p[t,t^{-1}]/ (\xi_m (t))$,
 and $f: X^3 \rightarrow X$ be
 defined by $f(x_1, x_2, x_3)=(x_1-x_2)^{a_1} (x_2 - x_3)^{a_2}$,
 where $a_i=p^{m_i}$,  for non-negative integers $m_i$, $i=1,2$.
 If $\xi_m (t)$ divides $1- t^{a_1+a_2}$, then
 \end{sloppypar}
 $$ \Phi_f(T(2,m)) = \{ \sqcup_{|X|^2} 
 ( - t \xi'_m (t) )^{a_1} s^{a_1 + a_2}
  \ | \  s \in X \}. $$
 \end{proposition}

\begin{example} {\rm
Using 
Proposition~\ref{3cocyinvprop} and the contribution formula preceding it, %
we obtain the following list of
calculations of $3$-cocycle invariants for $T(2,m)$ torus knots
carried out by a {\it Maple} program, for cocycles of the form
$f(x,y,z)=(x-y)(y-z)^p$.

\bigskip

\newpage 

\noindent
$\bullet$ $p=2$, $f(x,y,z)=(x-y)(y-z)^2$.
\begin{itemize}
\setlength{\itemsep}{-3pt}
\item [$\ast$]$m=3$: $16+48\ U^t$,
\item [$\ast$]$m=5$: $4096$,
\item[$\ast$] $m=7$: $262144$,
 \item[$\ast$] $m=9$: $4194304 +
12582912 \ U^{(t^4+t^7+1)}$,
\item [$\ast$]$m=11$: $1073741824$,
\item[$\ast$] $m=13$: $68719476736$,
 \item [ $\ast$]$m=15$:
$1099511627776 + 3298534883328\  U^{(t^{13}+t^{10}+t^7+t^4+t)}$.
\end{itemize}
$\bullet$ $p=3$, $f(x,y,z)=(x-y)(y-z)^3$.
\begin{itemize}
\setlength{\itemsep}{-3pt}
\item[$\ast$] $m=3$: $243+486\  U^{(2t+2)}$,
\item[$\ast$]$m=5$: $531441$,
\item[$\ast$] $m=7$:  $387420489$,
\item[$\ast$] $m=9$: $94143178827 +
188286357654\ U^{(2t^7+2t^6+t^4+t^3+2t+2)}$.
\end{itemize}
$\bullet$ $p=5$, $f(x,y,z)=(x-y)(y-z)^5$.
\begin{itemize}
\setlength{\itemsep}{-3pt}
\item[$\ast$]$m=3$:
$625+3750U^{(t+3)}+3750\ ( U^{(4t+2)}+U^{(3t+4)}+U^{(2t+1)} )$,
\item[$\ast$] $m=5$: $48828125 + 97656250\ (
U^{(4t^3+2t^2+2t+4)}+U^{(t^3+3t^2+3t+1)} ) $,
 \item[$\ast$] $m=7$: $3814697265625$.
\end{itemize}

} \end{example}

Similar calculations are made
 for twist knots using colorings
given in Example~\ref{colorex},
with polynomial
cocycles
$f(x,y,z)=(x-y)^{a_1} (y-z)^{a_2}$, to obtain the following.

\begin{proposition} \label{twist3cocyprop}
The 3-cocycle invariant of the twist knot $k(2n)$ with $2n+2$ crossings
as depicted in Fig.~\ref{twistknots}
is given by
\begin{displaymath}
\Phi_{f}(K)=
 \{ \sqcup_{|X|^2}   [-nt^{-a_1}+(1+n(1-t))^{a_1+a_2}] s^{2}
  \ | \  s \in X \}
\end{displaymath}
for $X=\Z_p[t, t^{-1}]/ (t - n (1-t)^2 )$ and
$f(x, y, z)=(x-y)^{a_1} (y-z)^{a_2}$.
\end{proposition}

\begin{example}\label{twistex}
{\rm
\begin{sloppypar}
The formula in Proposition~\ref{twist3cocyprop} is  input in
 {\it Maple} to obtain the following results.
\end{sloppypar}

\bigskip

\noindent
$\bullet$ $p=3$:
\begin{itemize}
\setlength{\itemsep}{-3pt}
\item[$\ast$] $n = 1$, $a_1=1$, $a_2=3$:
   $81 + 324\ (U^{(t + 2)}+ U^{(1 + 2 t)} ) $,
\item[$\ast$]$n =  1$, $a_1=3$, $a_2=1$:
$81+324\ (U^{(2t+2)}+U^{(t+1)})$,
 \item[$\ast$] $n=2$, $a_1=1$,
$a_2=3$: $243+486\ U^{(t+1)}$,
 \item[$\ast$] $n=2$, $a_1=1$, $a_2=3$:  $243+486\ U^{(2t+2)}$.
\end{itemize}
$\bullet$ $p=5$:
\begin{itemize}
\setlength{\itemsep}{-3pt}
\item[$\ast$] $n = 1$, $a_1=1$, $a_2=5$:
$3125+6250\ ( U^{(3t+3)}+U^{(2t+2)}) $,
 \item[$\ast$] $n=1$, $a_1=5$, $a_2=1$:
  $3125+6250\ (U^{(3t+3)}+U^{(2t+2)} )$,
   \item [$\ast$]$n=2$,
$a_1=1$, $a_2=5 $: $15625$,
 \item[$\ast$] $n=2$, $a_1=5$, $a_2=1$:
$15625$.
 \item[$\ast$]
$n=3$, $a_1=1$, $a_2=5$:
$625+3750\ (U^{(t)}+U^{(2t)}+U^{(3t)}+U^{(4t)})$,
\mbox{\hspace{1cm}}
 \item[$\ast$]
$n=3$, $a_1=5$, $a_2=1$:
$625+3750\ (U^{(t+1)}+U^{(2t+2)}+U^{(3t+3)}+U^{(4t+4)})$,
\item[$\ast$] $n=4$, $a_1=1$, $a_2=5$:
$625+3750\ (U^{(t+3)}+U^{(2t+1)}+U^{(3t+4)}+U^{(4t+2)})$.
\end{itemize}
$\bullet$ $p=7$:
\begin{itemize}
\setlength{\itemsep}{-3pt}
\item[$\ast$] $n=1$, $a_1=1$, $a_2=7$:
$2401+19208\ ( U^{(t+3)}+U^{(4t+5)}+U^{(2t+6)}
 +U^{(5t+1)}+U^{(6t+4)}+U^{(3t+2)}) $,
\item[$\ast$] $n=1$, $a_1=7$, $a_2=1$:
$2401+19208\ (U^{(t+1)}+U^{(2t+2)}+U^{(3t+3)}+U^{(4t+4)}+U^{(5t+5)}+U^{(6t+6)})$,
\item[$\ast$] $n=2$, $a_1=1$, $a_2=7$: $117649$,
 \item[$\ast$]
$n=2$, $a_1=7$, $a_2=1$ : $117649$,
 \item[$\ast$] $n=3$, $a_1=1$,
$a_2=7$:
$2401+19208\ (U^{(3t+4)}+U^{(5t+2)}+U^{(6t+1)}+U^{(t+6)}+U^{(2t+5)}+U^{(4t+3)} )$,
\item [$\ast$]$n=3$, $a_1=7$, $a_2=1$:
$2401+19208\ ( U^{(t+1)}+U^{(2t+2)}+U^{(3t+3)}+U^{(4t+4)}+U^{(5t+5)}+U^{(6t+6)} )$,
\item[$\ast$] $n=4$, $a_1=1$, $a_2=7$:
$2401+19208\ ( U^{(t+2)}+U^{(2t+4)}+U^{(3t+6)}+U^{(4t+1)}+U^{(5t+3)}+U^{6t+5)} )$,
\item[$\ast$] $n=4$, $a_1=7$, $a_2=1$:
$2401+19208\ ( U^{(t+1)}+U^{(2t+2)}+U^{(3t+3)}+U^{(4t+4)}+U^{(5t+5)}+U^{(6t+6)} )$,
\item[$\ast$] $n=5$, $a_1=1$, $a_2=7$:
$16807+33614\ ( U^{(3t+3)}+U^{(5t+5)}+U^{(6t+6)})$,
\item[$\ast$] $n=5$, $a_1=7$, $a_2=1$:
$16807+33614\ (U^{(t+1)}+U^{(2t+2)}+U^{(4t+4)}) $,
\item[$\ast$] $n=6$, $a_1=1$, $a_2=7$: $117649$,
 \item
[$\ast$]$n=6$, $a_1=7$, $a_2=1$: $117649$.
\end{itemize}
} \end{example}

\begin{remark}{\rm
Polynomial cocycles are further utilized in \cite{Chad}
for extensive computer calculations using {\it Maple}
and the knot table. The formulas given in
Propositions~\ref{2cocyinvprop} and \ref{3cocyinvprop},
however, enable one to compute for a larger quandles and knots.

Polynomial cocycles were also used in \cite{KAdissertation}
to evaluate cocycle invariants for twist spun knots using Alexander
quandles, using formulas in \cite{AS}.
Such computations are useful in detecting non-invertibility of
twist-spun knots.

} \end{remark}

\section{Non-triviality of Quandle Homology Groups}\label{homgpsec}

In this section we prove the non-triviality of homology groups of
some families of Alexander quandles.
The method
for dimensions $2$ and $3$ 
is to show that there is a coloring of knot diagrams that
evaluates non-trivially by a polynomial cocycle constructed
 in Lemma~\ref{polycocylemma}. This method has been used repeatedly
 since \cite{CJKLS}.
 For higher dimensions, we use algebraic machineries developed in \cite{NP},
 so that we give statements 
 and proofs separately.

\begin{proposition}\label{nontrivthm}
The following quandle homology groups $H^m_{\rm Q}(X;X)$
are non-trivial $(\neq 0): $
\begin{enumerate}
\setlength{\itemsep}{-3pt}
\item [{\rm (1)}]
 $X={\mathbb Z}_2[t, t^{-1}]/( \xi_{2^n + 1}(t))$ for any positive integer $n$,
 and for  $m=2,3$.
 \item [{\rm (2)}]
 $X={\mathbb Z}_p[t, t^{-1}]/(\xi_{(p^n + 1)/2} (t))$ for any odd prime $p$
 and  for any positive integer $n$, and for $m=2,3$.
 \item  [{\rm (3)}]
 $X=\mathbb{Z}_p[t,t^-1]/(t-n(1-t)^2)$, for $m=3$,
 and for:
 \begin{itemize}
 \setlength{\itemsep}{-3pt}
\item [$($a$)$] $p=3$, $n\equiv 1,2\pmod{3}$,
\item [$($b$)$] $p=5$, $n\equiv 1,3$, and $4$ $\pmod{5}$,
\item [$($c$)$] $p=7$, all $n$ except $n\equiv 2,6\pmod{7}$,
\item [$($d$)$] $p=11$, all $n$ except $n\equiv 1,2,6$, and $9\pmod{11}$,
\item [$($e$)$] $p=13$, all $n$ except $n\equiv 2,4,6,7$, and $12\pmod{13}$.
\end{itemize}
 \end{enumerate}
 \end{proposition}
{\it Proof.\/}  It is sufficient to show
that there is a coloring of a knot contributing a non-trivial value
to the quandle cocycle invariant~\cite{CJKLS}.

(1) $m=2,3$:  For $X={\mathbb Z}_2[t, t^{-1}]/ (\xi_{2^n + 1}(t))$,
let $f(x_1, x_2)=(x_1-x_2)^{2^n}x_2$ in Lemma~\ref{polycocylemma},
so that  $a_1=2^n$ and $a_2=1$. Note that
$1-t^{(a_1+a_2)}=(1-t)\xi_{2^n + 1} (t)$. Take $(1,0)\in X \times X$ as
a top color vector, which extends to a coloring of $T(2,m)$, where
$m=2^n + 1$ by Example~\ref{colorex}. Then by
Propositions~\ref{2cocyinvprop}, 
the contribution to the invariant is a multiple by a power of $t$ of
$\xi_{m}' (t)$, which is non-trivial in $X$, as the
degree of $\xi_{m}' (t)$ is less than that of $\xi_m (t)$.
For the $3$-cocycle case, choosing $a_1=1$ and $a_2=2$ gives
a contribution $\xi_{m}' (t)^{a_1}=\xi_{m}' (t)$ which is non-zero in $X$
and the same argument applies.

(2) In this case take $a_1=p^n$ and $a_2=1$ as before, then
$1-t^{(a_1+a_2)}=(1-t^{(p^n+1)/2})(1+t^{(p^n+1)/2})$. If $(p^n+1)/2$
is odd, then $\xi_{(p^n + 1)/2} (t)$ divides $(1+t^{(p^n+1)/2})$, and if
even, it divides $(1-t^{(p^n+1)/2})$, hence the result follows by
the same argument.



(3) From Proposition~\ref{twist3cocyprop},
the invariant is non-trivial if
 $[-nt^{-a_1}+(1+n(1-t))^{a_1+a_2}]$ that appear
in the formula is non-zero.
The choice $a_1=1$ and $a_2=p$ gives
non-zero values for the cases listed in the statement,
among primes $2<p\leq 29$, and for $0<n<p$.
$\hfill$ $\Box$

\bigskip

To prove non-triviality for higher dimensions we use the following lemma
from \cite{NP}.
For a quandle $(X, *)$ and a positive integer $m$,
we use the notation $*_a (x)=x*a$ and
$(*_a)^m (x)=(\cdots ( x *a) *a) *\cdots )*a$ where the operation is performed
$m$ times for a positive integer $m$ and for $x, a \in X$.
For an $n$-chain $c=(x_1 , \ldots, x_n)$, $x_i \in X$, $i=1, \ldots, n$
for a positive integer $n$, the notation $*_a (c)=c*a=(x_1 *a, \ldots, x_n*a)$
is used. The map is extended to the chain groups linearly.
The map $h_a: C_n^{\rm Q}(X) \rightarrow C_{n+1}^{\rm Q}(X) $
was defined in \cite{NP} by linearly extending
$$h_a(c)=h_a( (x_1 , \ldots, x_n) )=(x_1 , \ldots, x_n, a)=(c,a). $$
Then $h'_a: C_n^{\rm Q}(X) \rightarrow C_{n+1}^{\rm Q}(X) $ is defined by
 $h'_a=h_a+*_a h_a+\cdots+ (*_a)^m  h_a$ for $a \in X$.

Let $s=s(y_0, y_1)=\sum_{i=1}^{k-1} (y_{i-1}, y_{i})$ for
$y_i \in X$, $i=0, \ldots, k-1$. Define \cite{NP}
 $h_s:C_n^{\rm Q}(X) \rightarrow C_{n+2}^{\rm Q}(X) $ by
 linearly extending
 $$h_s(c)=(c, s)=\sum_{i=1}^{k-1} (x_1 , \ldots, x_n, y_{i-1}, y_{i} ) . $$

\begin{lemma}[\cite{NP}] \label{chainmaplemma}
{\rm (i) }
Let $X$ be a quandle such that there is an element $a \in X$
 satisfying the condition
$(*_a)^m (x) =x$ for any $x \in X$, then $h'_a$ is a chain map.

\noindent
{\rm (ii)} Let $X$ be a quandle such that there is a sequence of elements
$(y_0, \ldots, y_k)$ satisfying the condition
$y_{i+1}=y_{i-1}*y_i$ for $i=1, \ldots, k-1$,
$y_0=y_{k-2}*y_{k-1}$ and $y_1=y_{k-1}*y_{0}$.
Then $h_s$ is a chain map.
\end{lemma}

\begin{theorem}
For any positive integer $n$,
the quandle $X=\mathbb{Z}_2[t,t^{-1}]/ ( \xi_{m} (t) )$,
where $m=2^n+1$,
has non-trivial $4$-dimensional cohomology:
$H^4(X;X)\neq 0$.
\end{theorem}
{\it Proof.\/}
We exhibit a cycle $C$ and a cocycle $f$ such that $f(C)\neq 0$ to prove
non-triviality. 
Although negative signs are irrelevant in $\Z_2$, we often  leave
them to indicate computational processes below. 
 We
construct a $4$-cycle using a $3$-cycle and the map defined above.
From a shadow coloring of a $(2,m)$-torus knot with the quandle
$X$
such that the top color vector is $(0,1)$ and the left-most region 
is colored by $0$, 
we have a 3-cycle $C_3=\sum_{k=0}^{m-1}(0,\xi_k, \xi_{k+1})$,
and define $h_0(C_3)=C_4$  and $h'_0(C_3)=C'_4$. By
Lemma~\ref{chainmaplemma}, $C'_4$ is a $4$-cycle. We use the
$4$-cocycle
$$f(x_1,x_2,x_3,x_4)=(x_1-x_2)^{a_1}(x_2-x_3)^{a_2}(x_3-x_4)^{a_3},
\quad a_i=2^{n_i}, \quad  i=1,2,3. $$
Using $$(x_i*b - x_{i+1}*b)^{a_i} =( (tx_i + (1-t)b ) - ( tx_{i+1} + (1-t)b))^{a_i}=
t^{a_i} (x_i - x_{i+1}), $$
we have
 $$f(C'_{4})=(1+t^{a_1+a_2+a_3}+\cdots +
t^{(m-1)(a_1+a_2+a_3)})
f(C_4), $$
where $f(C_4)=\sum_{k=0}^{m-1}(0-\xi_k)^{a_1}(\xi_k-\xi_{k+1})^{a_2}\xi_{k+1}^{a_3}$.
Suppose that $1-t^{a_1+a_2+a_3}=0$, then
\begin{eqnarray*}
f(C'_4)&=& m\sum_{k=0}^{m-1}(-t)^{ka_2}\xi_k^{a_1}(1-t\xi_k)^{a_3}\\
&=&
m\sum_{k=0}^{m-1}(-t)^{ka_2}\xi_{k}^{a_1}-t^{a_3}\sum_{k=0}^{m-1}(-t)^{ka_2}\xi_k^{a_1+a_3} .
\end{eqnarray*}
Set  $a_1=a_2=2^{n-1}$, $a_3=1$ and
$m=2^n+1$, then $ 1-t^{a_1+a_2+a_3}=1-t^{2^n+1}=0$.
Then
\begin{eqnarray*}
f(C'_4)&=&(\sum_{k=0}^{m-1}(-t)^k\xi_k)^{a_1}-t^{a_3}\sum_{k=0}^{m-1}(-t)^{ka_2}\xi_k^{a_1+a_3}.
\end{eqnarray*}
By induction we see  $\displaystyle \sum_{k=0}^{m-1}(-t)^k\xi_k=\frac{\xi_m
\xi_{m+1}}{1-t}$ (shown in  \cite{KAdissertation}).
 Since $m$ is odd, $\xi_{m+1}$ in the RHS is divisible by
 $(1-t)$, and  this sum is $0$ in $X$. Then in
$X$
\begin{eqnarray*}
f(C'_4)&=&-t^{a_3}\sum_{k=0}^{m-1}(-t)^{ka_2}\xi_k^{a_1+a_3}.
\end{eqnarray*}
We now compute the sum $ \displaystyle
\sum_{k=0}^{m-1}(-t)^{ka_2}\xi_k^{a_1+a_3}$. Note that $1+t$ is invertible
in $X=\mathbb{Z}_2[t,t^{-1}]/ (\xi_m (t)) $,
since
$1+t+t^2+\cdots + t^{m-1}=0$ implies
$(1+t)(t+t^3+\cdots+ t^{m-2})=1$.
Hence $(1+t^{a_1})=(1+t)^{a_1}$ and $(1+t^{a_3})$ are invertible.
Then we compute
\begin{eqnarray*}
(1+t^{a_1})(1+t^{a_3})\sum_{k=0}^{m-1}(-t)^{ka_2}\xi_k^{a_1+a_3}&=&
\sum_{k=0}^{m-1}(-t)^{ka_2}(1-(-t)^{ka_1})(1-(-t)^{ka_3})\\
&=&\sum_{k=0}^{m-1}(-t)^{k(a_1+a_2+a_3)}=m
\end{eqnarray*}
since $\sum_{k=0}^{m-1}(-t)^{ka_2}=(\xi_m)^{a_2}=0$,
 $\sum_{k=0}^{m-1}(-t)^{ka_3}=(\xi_m)^{a_3}=0$,
 and 
 $(-t)^{ k (a_2+ a_3)}=(-t)^{k (- a_1)}$.
Thus
$ f(C'_4)=-mt^{a_3}(1+t^{a_1})^{-1}(1+t^{a_3})^{-1} \neq 0$ in $X.$
$\hfill$ $\Box$

\begin{theorem} \label{mainthm}
For any positive integer $n$ and $r>1$, the quandle
$X=\mathbb{Z}_p[t,t^{-1}]/ ( \xi_{m} (t) )$ has
non-trivial $(2r+1)$-dimensional cohomology, $H^{2r+1}(X;X)\neq 0$,
if:

\vspace{-10pt}

\begin{itemize}
\setlength{\itemsep}{-3pt}
\item[{\rm (i)}] $p=2$ and  $m=2^n+1$, or, 
\item[{\rm (ii)}]
 $p$ is an odd prime and $m=(p^n+1)/2$.
 \end{itemize}
\end{theorem} {\it Proof.\/}
 We construct a $(2r+1)$-cycle using the $(2r-1)$-cycle
 $C_{2r-1}$  and the chain map $h_s$ defined above, where
  $s=\sum_{k=0}^{m-1}(\xi_k, \xi_{k+1})$.
 Let $C_{2r+1}=h_s(C_{2r-1})$.
 To compute the cocycle values, denote
 $C_{2r-1} $ by a formal sum
 $\sum (x_1, \ldots, x_{2r-1})$ and
 $C_{2r+1}$ by $\sum_{k=0}^{m-1}   \sum (x_1, \ldots, x_{2r-1}, \xi_k, \xi_{k+1})$.
By  Lemma~\ref{chainmaplemma},
 $C_{2r+1}$ is a $(2r+1)$-cycle. We use the
$(2r+1)$-cocycle
$$f_{2r+1}(x_1, \ldots, x_{2r+1})=(x_1-x_2)^{a_1}\cdots (x_{2r}-x_{2r+1})^{a_{2r}},
$$
for $k=1, \ldots, r$. Then one computes
\begin{eqnarray*}
f(C_{2r+1})&=& \sum_{j=0}^{m-1}\sum
 f(x_1, \ldots, x_{2r-1}, \xi_j, \xi_{j+1})\\ 
&=&
\sum_{j=0}^{m-1}\sum(x_1- x_2)^{a_1}\cdots (x_{2r-2}-x_{2r-1})^{a_{2r-1}}
(x_{2r-1} - \xi_{j})^{a_{2r-1}}   (\xi_j-\xi_{j+1})^{a_{2r}}\\
&=&
\sum_{j=0}^{m-1}\sum(x_1- x_2)^{a_1}\cdots (x_{2r-2}-x_{2r-1})^{a_{2r-1}}
(x_{2r-1})^{a_{2r-1}}   (\xi_j-\xi_{j+1})^{a_{2r}} \\
& -&\sum_{j=0}^{m-1}\sum(x_1- x_2)^{a_1}\cdots (x_{2r-2}-x_{2r-1})^{a_{2r-1}}
( \xi_{j})^{a_{2r-1}}   (\xi_j-\xi_{j+1})^{a_{2r}} . 
\end{eqnarray*}
The first term vanishes because
$\sum_{j=0}^{m-1}  (\xi_j-\xi_{j+1})^{a_{2r}}=\xi_m^{a_{2r}}  =0$ in $X$.
If $a_{2r-1}$ is odd, then 
the second term is
$$\sum(x_1- x_2)^{a_1}\cdots (x_{2r-2}-x_{2r-1})^{a_{2r-1}}
\sum_{j=0}^{m-1}(0 -  \xi_{j})^{a_{2r-1}}
(\xi_j-\xi_{j+1})^{a_{2r}}.$$ Now we prove the theorem by
induction by proving that $f_{2r+1}(C_{2r+1})$ is invertible in
$X$ under the assumption  $f_{2r-1}(C_{2r-1})$ is invertible. 

For the case {\rm (i)}, let $a_{2k-1}=1$ and $a_{2k}=2^n$. 
For
$r=1$, $f_3 (C_3)=\sum_{j=0}^{m-1}(0 -  \xi_{j})
(\xi_j-\xi_{j+1})^{2^n}= t\xi_m'$, and $t (1+t)\xi'_m(t)=1$ in $X$. 
Assume $f_{2r-1}(C_{2r-1})= \sum_{j=0}^{m-1}\sum(x_1-
x_2)^{a_1}\cdots (x_{2r-2}-x_{2r-1})^{a_{2r-1}}$ is invertible in
$X$. Then the above second term is  $f_{2r-1}(C_{2r-1}) f_3 (C_3)$,
which is invertible in $X$  
and this case is proved. 

For the case {\rm (ii)}, let $a_{2k-1}=1$ and $a_{2k}=p^n$. 
 For
$r=1$, $f_3 (C_3)=\sum_{j=0}^{m-1}(0 -  \xi_{j})
(\xi_j-\xi_{j+1})^{p^n}= t\xi_m'$.
Differentiating both sides of $(1+t) \xi_m= 1- (-t)^m$, 
we obtain $(1+t)\xi'_m=-\xi_m-m(-1)^m t^{m-1}=-mt^{-1}$ in $X$.
By the assumption $m$ is invertible in $\mathbb{Z}_p$, 
hence $t\xi'_m$ is invertible in $X$. 
Assume $f_{2r-1}(C_{2r-1})= \sum_{j=0}^{m-1}\sum(x_1-
x_2)^{a_1}\cdots (x_{2r-2}-x_{2r-1})^{a_{2r-1}}$ is invertible in
$X$. Then the above  term   $f_{2r-1}(C_{2r-1}) f_3 (C_3)$
is invertible in $X$ 
and the theorem is proved.
$\hfill$ $\Box$
\bigskip

\begin{sloppypar}
\noindent 
{\bf Acknowledgments.\/} We are thankful to J.\,S.\, Carter, M.\, Elhamdadi, 
M.\, Niebrzydowski, and J.\,H.\, Przytycki for
continuous valuable conversations.
\end{sloppypar}

\end{document}